\documentclass[12pt]{article}
\newcommand{\R}{\mbox{I$\!$R}}

\newcommand{\bk}{\\[0.03in] \hspace*{0.5in} }
\newcommand{\btd}{\bigtriangledown}

\newcommand{\mfor}{\ \ \ \ {\mbox{for}} \ \ }

\begin{document}

\begin{center} {\LARGE   {\bf Blow-up Solutions of }}
\medskip \smallskip \\ {\LARGE {\bf Nonlinear Elliptic Equations in ${\R}^n$}} 
\medskip \smallskip \\ {\LARGE {\bf with Critical Exponent}} 
\medskip \medskip  \medskip  \smallskip \\ 
{\large {Man Chun {\Large L}EUNG}}\\
\smallskip
{\large {Department of Mathematics,}}\\
{\large {National University of Singapore,}}\\
{\large {Singapore 117543}}\\
{\large {Republic of Singapore}}\\
{\tt matlmc@math.nus.edu.sg}
\end{center}
\vspace{0.4in}
\begin{abstract}
\noindent For an integer $n \ge 3$ and any positive number $\varepsilon$ we establish the existence of smooth functions $K$ on $\,{\R}^n \setminus \{ 0 \}\,$ with $\,|K - 1| \le \varepsilon,\,$ such that 
the equation 
$$
\Delta u + n (n - 2) K u^{{n + 2}\over {n - 2}} = 0
\ \ \ \ {\mbox{in}} \ \ {\R}^n \setminus \{ 0 \} 
$$
 has a smooth positive solution which blows up at the origin (i.e., $u$ does not have slow decay near the origin).  
Furthermore, we show that in some cases $K$ can be extended as a Lipschitz function on ${\R}^n.$ These provide counter-examples to a conjecture of C.-S. Lin  when $n > 4\,,\,$ and Taliaferro's  conjecture.

\end{abstract}

\vspace*{1in}

KEY WORDS:   nonlinear differential equation, blow-up, positive solutions, decay estimates.\\[0.1in]
2000 AMS MS CLASSIFICATIONS: Primary 35J60\,; Secondary 53C21.

\pagebreak

{\bf \large {\bf 1. \ \ Introduction}}

\vspace{0.2in}

Given any positive number $\varepsilon > 0\,,$ we construct blow-up solutions of the equation 
$$
\Delta u + n (n - 2) K u^{{n + 2}\over {n - 2}} = 0 \ \ \ \ {\mbox{in}} \ \ {\R}^n \setminus \{ 0 \}\,, \leqno (1.1)
$$
with $|K - 1| \le \varepsilon$ in ${\R}^n \setminus \{ 0 \}\,.$ Here $K$ is a smooth function on ${\R}^n \setminus \{ 0 \}\,,$ and $n \ge 3$ is an integer. A positive smooth solution $u$ of equation (1.1) is known as a blow-up solution if $u$ does not have slow decay near the origin, i.e., there does {\it not} exist a positive constant $C$ such that the slow decay estimate
$$
u (x) \le C |x|^{{(2 - n)}\over 2} \leqno (1.2)
$$
is satisfied for all $x \in {\R}^n \setminus \{ 0 \}$ with $|x|$ close to zero. 
The question on slow decay is a crucial first step in understanding equation (1.1) regarding   such basic questions as existence of solutions, asymptotic behavior and properties of the moduli space (cf. \cite{Chen-Lin.1}, 
\cite{Chen-Lin.4}, \cite{Gromov} and \cite{M-P-U}).\bk
Equation (1.1) is recognized as a conformal scalar curvature equation, as it provides via $\,4n (n - 1)K$ the scalar curvature of the conformal metric $g_c := \displaystyle{u^{4\over {n - 2}} \, \delta_{ij}\,.}$ Many authors study the equation from the view points of nonlinear partial differential equations, prescribing curvature, parabolic flow and asymptotic geometry. Surprising and deep features of the equation are revealed (see, for instances, \cite{Ambrosetti_et_al}, \cite{Aubin.1}, \cite{Aubin_Bahri}, \cite{Aubin_Bismuth}, \cite{Chen-Lin.2}, \cite{Chen-Lin.3},   \cite{ChenWX-Li-Ann}, \cite{Kato}, \cite{Kazdan}, \cite{K-M-P-S}, \cite{Lee-Parker}, \cite{Leung.1}, \cite{Leung.5}, \cite{Li}, \cite{Lin.1}, \cite{Loewner-Nirenberg}, \cite{M-P},  \cite{Schoen},  \cite{Yanagida_Yotsutani}, \cite{Zhang} and the references within). \bk
One of the early important developments on the equation is made by Gidas, Ni and Nirenberg (\cite{Gidas-Ni-Nirenberg.1}, \cite{Gidas-Ni-Nirenberg.2}; cf. \cite{Caffarelli-Gidas-Spruck}), who show that if equation (1.1) is satisfied on the whole ${\R}^n$ with $K \equiv 1\,,$ then 
$$
u (x) = \left( {\lambda\over {\lambda^2 + |x - \xi|^2}} \right)^{{n - 2}\over 2} \mfor x \in {\R}^n, \leqno (1.3)
$$
where $\lambda$ is a positive constant and $\xi$ a fixed point in ${\R}^n.$ We call solutions of the form (1.3) {\it standard bubbles} and $\xi$ the {\it center} of the standard bubble, partly because the conformal metric $g_c$ is isometric to a sphere in this case.\bk 
Later, by using the moving plane method propounded by the above authors, Caffarelli, Gidas, and Spruck \cite{Caffarelli-Gidas-Spruck} consider the case when $K \equiv 1$ in a non-empty open neighborhood of the origin and find out that either $u$ has a removable singularity at the origin, or it has slow decay (1.2) and  precise asymptote near the origin as described by the Delaunay-Fowler-type  solutions (cf. section 2). Inspired by these exquisite results, people begin to look at what may happen when $K$ is close to, but not identically equal to, one.\bk 
In an enlightening article by Korevaar, Mazzeo, Pacard and  Schoen \cite{K-M-P-S}, the problems are reconsidered by introducing subtle geometric insights which lead to a second order term in the asymptotic expansion. Included in that paper is the observation 
that if $u$ fails to have slow decay (1.2), then, roughly speaking, $u$ is close to the standard bubble with very small $\lambda$ in a neighborhood of a local maximum near the origin. This is commonly known as the blow-up process.\bk
Blow-up solutions are first constructed by Taliaferro \cite{Taliaferro.1} with $K$ being non-constant and bounded between two positive constants $a^2$ and $b^2 \ (b > a > 0)\,.$ The method is by superimposing standard bubbles in a controlled way, which is applied by the author in \cite{Leung.6} for further global examples. Taliaferro also asks whether when $\displaystyle{ b/a < 2^{2\over {n - 2}}}\,,$  then $u$ has slow decay (\cite{Taliaferro.1}; the number $\displaystyle{ 2^{2\over {n - 2}} }$ appears in bounding $K$ when two standard bubbles are added). Our construction gives a negative answer to the question. Nevertheless, fruitful studies that contribute to deeper understanding on equation (1.1) are stimulated by the question (cf. 
\cite{Cheung-Leung.1}, \cite{Leung.6}, \cite{Leung.9} and \cite{Taliaferro.2}).\bk
When $K$ is non-constant but has certain order of expansion at the origin, Chen and Li in a series of masterful papers \cite{Chen-Lin.2} \cite{Chen-Lin.3} \cite{Lin.1} investigate the question of slow decay, among other important issues. Li conjectures that if $K$ is positive and H\"older-continuous on ${\bar B}_2$, then $u$ has slow decay \cite{Lin.1}. Here ${\bar B}_2$ is the closed ball in ${\R}^n$ with center at the origin and radius 2.  We also give a negative answer to the conjecture when $n > 4$ by finding a smooth positive function $K$ on ${\R}^n \setminus \{ 0\}$ which can be extended to be a Lipschitz function on ${\R}^n$ by setting $K (0) =1\,,$ such that equation (1.1) has a blow-up solution. At the present moment, the conjecture remains open in 3 and 4 dimensions.\bk
It appears that both the standard bubbles and the Delaunay-Fowler-type  solutions are too rigid to be `deformed' substantially without altering $K$ a lot. Rather, they follow the rules of conformal transformations such as translations and reflections (Kelvin transforms). We first make use of the periodicity of Delaunay-Fowler-type  solutions to modify them with the standard bubbles (section 2). As each Delaunay-Fowler-type  solution has slow decay, we seek to offset the `centers'. This is based on the observation that a standard bubble is symmetric with respect to specific Kelvin transforms (sections 4 and 5). The construction is then by putting the offset and modified Delaunay-Fowler-type  solutions on a collection of small open balls with the origin as a limit point. The resulting $K$ is 
equal to one except on a set of finite  measure, which can be made arbitrarily small. We also 
find out that in some cases $K$ can be extended as a Lipschitz function on ${\R}^n$ (section 7).    

\vspace*{0.3in}
{\bf Acknowledgment.} \ \ The author is pleased to express his gratitude to Rafe Mazzeo for delightful discussions and insights into the equations, which encouraged the author's endeavour; and to C.-S. Lin for sharing his view on many problems and hospitality during the author's visit to the National Center for Theoretical Sciences (Taiwan).

\vspace{0.3in}

{\bf \large {\bf 2. \ \ Modified Delaunay-Fowler-type solutions}}

\vspace{0.2in}

Consider equation (1.1) in the cylindrical coordinates $(t, \,\,\theta)\,,$ where 
$$
t = - \ln |x| \ \ \ \ {\mbox{and}} \ \ \ \ \theta = x/|x| \mfor x \in {\R}^n \setminus \{ 0 \}\,.
$$ 
The function 
$$
v (t, \theta) := |x|^{{n - 2}\over 2} u (x) \mfor t \in \R \ \ {\mbox{and}} \ \ \theta \in S^{n - 1}  \leqno (2.1)
$$
satisfies the equation 
$$
{{\partial^2 v}\over {\partial t^2}} + \Delta_\theta \,v - {{(n - 2)^2}\over 4} \,v 
+  n (n - 2) K_{cy} \,v^{{n + 2}\over {n - 2}} = 0 \ \ \ \ {\mbox{in}} \ \ \R \times S^{n - 1}.$$
Here $\Delta_\theta$ is the Laplacian on the standard  unit sphere in ${\R}^n$ and $K_{cy} (t, \theta) := K (x)\,,$ with $|x| = e^{-t}$ and $x/|x| = \theta$. \bk 
We are interested in {\it radial} solutions of the above equation, especially when $K_{cy} \equiv 1$:
$$
v'' - {{(n - 2)^2}\over 4} \,v 
+  n (n - 2)   \,v^{{n + 2}\over {n - 2}} = 0 \ \ \ \ {\mbox{in}} \ \ \R\,.
\leqno (2.2)
$$ 
Positive smooth solutions of (2.2) are known as Delaunay-Fowler-type solutions. By a re-parametrization, we standardize the solutions so that 
$$
v (0) = {\mathop{\max}\limits_{t \in \R}} \ v (t)\,. \leqno (2.3)
$$ 
Delaunay-Fowler-type   solutions can be indexed by the neck-size which is equal to $\eta := {\mathop{\inf}\limits_{t \in \R}} \ v (t)\,.$ For $\eta  > 0$, $v$ is periodic \cite{M-P}. The  relation between $\eta$ and the period $T$ is described in the following (an equivalent form is derived in \cite{M-P}; see also \cite{Leung.8}).\\[0.2in]
{\bf Proposition 2.4.} \ \ {\it There exists a positive constant $C$, independent on $\eta$ and $T\,,$ such that }
$$
C^{-1} e^{-  \left( {{n - 2}\over 4} \right) T}
 \le \eta \le C e^{-  \left( {{n - 2}\over 4} \right) T}
 \ \ \ \ for \ \ T \gg 1\,. \leqno (2.5)
$$

\vspace*{0.2in}

\hspace*{0.5in}We find it more convenient in our case to index the Delaunay-Fowler-type solutions by the period $T\,.$ Thus for $T \gg 1\,,$  
denote by $v_T$ the Delaunay-Fowler-type   solution with period $T$, neck-size $\eta_T\,,$ and with standardization (2.3). 
Given any positive number $D$, it is known that $v_T$ converges uniformly in $C^2$-norm on $[-2D\,, \ 2D]$ to the solution 
$$
v_s (t) : = \displaystyle{(2 \cosh \,t)^{(2 - n)/2}} \leqno (2.6)
$$
of (2.2) (see \cite{K-M-P-S} and \cite{M-P}; cf. also \cite{Gilbarg-Trudinger} and \cite{Leung.8}). Because of its simple and explicit form, we seek to replace part of $v_T$ by $v_s\,.$\bk
Given a (small) positive number $\varepsilon$ and a (large) number $D$, choose a large number $T$ so that $T \gg 2D$ and $v_T$ is close to the function $\displaystyle{(2 \cosh \,t)^{(2 - n)/2}}$ on $[-2D\,, \ 2D]\,.$ We start to perform a cut-and-glue-in process on $v_T$ on $(D,\  2D]$ so that the modified function ${\breve{v}}_T$ remains smooth and positive and 
\[{\breve{v}}_T (t) \ = \ \left\{ \begin{array}   
{r@{\quad \mbox{for}\quad}l}
 (2 \cosh \,t)^{(2 - n)/2} & \ \  t \le D\,,\\  
 v_T (t)  \ \ \ \ \ \ \ & \ \  t \ge 2D\,. 
\end{array} \right.  \]
We refer to the appendix for details on the construction of ${\breve{v}}_T\,.$  
By choosing $T$ large enough, we may also assume that ${\breve{v}}_T$ satisfies the equation 
$$
{\breve{v}}_T'' - {{(n - 2)^2}\over 4} {\breve{v}}_T + n (n - 2) {\breve K} \,{\breve{v}}_T^{{n + 2}\over {n - 2}} = 0 \ \ \ \ {\mbox{in}} \ \ \R\, \leqno (2.7)
$$
with $|\breve K - 1| \le \varepsilon$ in $\R$ ($\breve K (t) = 1$ for $t \not\in [D, \ 2D]$). Following closely the argument in \cite{Leung.8}, we show in the appendix that 
$$
|\breve K (s) - \breve K  (t) | \le C'_D\, \eta_T^2\, |s - t| \mfor s, \ t \in \R\,, \leqno (2.8)
$$
and for $T$ large enough ($\eta_T$ is the neck-size). Here $C'_D$ is a positive constant that depends on $D$ and $n$ only. \bk
From the periodicity of $v_T$, $v_T (t - T) = v_T (t)$ for $t \in \R\,.$ It follows that $v_T (t)$ is close to the function $\displaystyle{\left(2 \cosh \, (t - T)\right)^{(2 - n)/2}}$ on $[T -2D\,, \ T + 2D]\,.$
We perform the cut-and-glue-in process on $[T- 2D, \ T - D)$  so that the modified function, which is also denoted by ${\breve{v}}_T$ for the sake of brevity, is smooth and positive, and  
\[{\breve{v}}_T (t) \ = \ \left\{ \begin{array}   
{r@{\quad \mbox{for}\quad}l}
 v_s (t) \ \ \ \ \ \ \  \ \ \ \ \,\,\,\,& \ \  t \le D\,,\\ 
v_s (t - T)  \ \ \ \ \ \ \ & \ \  t \ge T - D\,, 
\end{array} \right.  \]
and (2.7) remains satisfied with (a different smooth function which we still denote by $\breve K\,,\,$ and) $|\breve K - 1| \le \varepsilon$ in $\R\,.$\bk
For later applications, the above construction suffices. The periodicity of $v_T$ allows us to carry on     
the cut-and-glue-in process $m$ times, and  obtain the modified Delaunay-Fowler-type solution, which is continued to be denoted by $\breve v_T\,$  ($m$ is an integer bigger than one). It has the following properties: 
$$
{\breve{v}}_T (t) = v_s (t) \mfor t \le  D\,, \leqno (2.9)
$$
$$
{\breve{v}}_T (t) = v_s (t - (k - 1)T) \mfor t \in [(k - 1)T - D\,, \ (k - 1)T + D]\,,  
\leqno (2.10)
$$
$$
{\breve{v}}_T (t) = v_s (t - (m - 1)T) \mfor t \ge (m -1)T - D\,, \leqno (2.11)
$$
where $k = 1, \ 2,\cdot \cdot \cdot , \ m-1\,\,,$ and
$$
{\breve{v}}_T'' - {{(n - 2)^2}\over 4}\, {\breve{v}}_T + n (n - 2) {\breve K} \,{\breve{v}}_T^{{n + 2}\over {n - 2}} = 0 \ \ \ \ {\mbox{in}} \ \ \R\,. \leqno (2.12)
$$
Here we continue to use $\breve K$ to denote the smooth function on $\R$ in (2.12) which satisfies  $|\breve K - 1| \le \varepsilon$ in $\R\,.$  By using the triangle inequality and  the fact that $\breve K = 1$ outside a compact subset, it can be seen that (2.8) remains valid.\bk
In order to find out how big the function may become in Euclidean coordinates, we reverse the process of cylindrical transform (2.1). Thus we seek to express parts of the modified Delaunay-Fowler-type   solution $\breve v_T$ in Euclidean coordinates. Consider the inverse  transformation of (2.1) 
$$
u (x) = |x|^{- {{n - 2}\over 2}} v (t, \theta) \mfor x\not= 0\,. \leqno (2.13)
$$
(Here, as above, $|x| = e^{-t}$ and $\displaystyle{\theta = {x\over {|x|}}}$ for $x \in {\R}^n \setminus \{ 0\}\,.$)\, Under  (2.13), the function $v_s (t) = (2 \cosh \,t)^{(2 - n)/2}$ corresponds to the {\it canonical standard bubble}   
$$
u_s (x) := \left( {1\over {1 + |x|^2}} \right)^{{n - 2}\over 2} \mfor x \in {\R}^n \setminus \{ 0\}, \leqno (2.14)
$$
which has a removable singularity at the origin. From \cite{Leung.7} we have the following result (mindful of the change on the sign of $t$).\\[0.2in]
{\bf Lemma 2.15} \ \ {\it Let 
$$
u_1 (x) = \left( {\lambda_1\over { \lambda_1^2 + |x|^2}} \right)^{{n - 2}\over 2}
\ \ \ \ {\mbox{and}} \ \ \ \ 
u_2 (x) = \left( {\lambda_2\over { \lambda_2^2 + |x|^2}} \right)^{{n - 2}\over 2} \leqno (2.16)  
$$
for $x \in \R\,,$ and $v_1$ and $v_2$ be the transforms of $u_1$ and $u_2$ according to (2.13), respectively. Then $v_1 (t) = v_2 (t - \overline t)$ for $t \in \R$, \,with\, 
$\overline t = \ln \,(\lambda_2/\lambda_1)\,.$ Here $\lambda_1$ and $\lambda_2$ are positive constants.}\\[0.2in]
\hspace*{0.5in}The proof is based on pointwise calculation (see \cite{Leung.7}), and the relation in lemma 2.15 can be used locally. Let ${\breve u}_T$ be the transform of ${\breve v}_T$ under (2.13). From (2.9) and (2.14), we have 
$$
{\breve u}_T (x) = u_s (x) =  \left( {1\over {1 + |x|^2}} \right)^{{n - 2}\over 2}\mfor |x| \ge e^{-D} \ \ \  \ (D \gg 1)\,, \leqno (2.17)
$$
which indicates that ${\breve u}_T$ is the same as the canonical standard bubble $u_s$ except on a small ball.  
As
${\breve v}_T (t) = v_s (t - T)$ for $t \in [T -D\,, \ T + D]\,,$ 
if we take $u_2 = u_s$ in lemma 2.15, we see that, under the transformation (2.13), ${\breve v}_T$ is expressed as 
$$
\left( {e^{-T}\over { e^{-2T} + |x|^2}} \right)^{{n - 2}\over 2} \mfor e^{-(T+D)} \le |x| \le e^{-T + D}\,.
$$
That is,
$$
{\breve u}_T (x) = \left( {e^{-T}\over { e^{-2T} + |x|^2}} \right)^{{n - 2}\over 2} \mfor e^{-(T+D)} \le |x| \le e^{-T + D}\,.
$$
In particular, take a point $x_T \in {\R}^n$ with $|x_T| = e^{-T}\,,$ we have
$$
{\breve u}_T (x_T) = \left( {{e^T}\over 2} \right)^{{n - 2}\over 2}\,, \leqno (2.18)
$$
which can be made arbitrarily large by choosing large enough $T\,.$\bk 
Similarly, when $m \ge 3\,,$ we have (2.17) and 
$$
{\breve u}_T (x) = \left( {\lambda_k\over { \lambda_k^2 + |x|^2}} \right)^{{n - 2}\over 2} \mfor e^{-(k - 1)T-D} \le |x| \le e^{-(k - 1)T + D}\,,
$$
together with 
$$
{\breve u}_T (x) = \left( {\lambda_m\over { \lambda_m^2 + |x|^2}} \right)^{{n - 2}\over 2} \mfor   0 < |x| \le e^{-(m - 1)T + D}\,. \leqno (2.19)
$$
Here $\lambda_k := e^{-(k - 1)T},\, \ \lambda_m := e^{-(m - 1)T},\,\,$ and $k = 2,\cdot \cdot \cdot ,\, m - 1\,.$  
It follows from (2.19) that ${\breve u}_T$ has a removable singularity at the origin. Thus by  choosing either $m$ or $T$ large (cf. (2.18)), ${\breve u}_T (0)$ can be as big as we like.\bk
 When $m = 2$, 
\begin{eqnarray*}
{\breve u}_T (x) & = &  \left( {1\over {1 + |x|^2}} \right)^{{n - 2}\over 2}  
 \mfor |x| \ge e^{-D},\\
{\breve u}_T (x) & = &  \left( {{e^{-T}}\over {e^{-2T} + |x|^2}} \right)^{{n - 2}\over 2}  
 \mfor 0 <  |x| \le e^{-T + D}, 
\end{eqnarray*}
and ${\breve u}_T (0) := e^{(2-n) \,T/2}\,.$ Provided that $T$ and $D$ are {\it large}, the graph of ${\breve u}_T$ is alike a \,(very tall and thin) needle stuck vertically on a (almost round) piece of  cotton. Geometrically, the conformal metric ${\breve u}_T^{4/(n - 2)} \delta_{ij}$ can be thought of as two spheres glued together near a common point, with the metric perturbed on a small neighborhood of the common point.\bk
For $k = 1, \ 2,\cdot \cdot \cdot , \,m\,,$ the standard bubbles 
$$
\left( {\lambda_k\over { \lambda_k^2 + |x|^2}} \right)^{{n - 2}\over 2} \mfor x \in {\R}^n \ \ \ \left(\lambda_k := e^{-(k - 1)T}\right),  
$$
are all rotationally symmetric above the origin. 
Overall, as noted in \cite{K-M-P-S}, the metrics $\,\displaystyle{v_T^{4\over{n - 2}} (dt^2 + d\theta^2)}\,$ converge as $T \to \infty$ to a bead of spheres of same radius which are arranged along a {\it fixed} axis. As each $v_T$ has slow decay (cf. (2.13)), in order to produce fast growing solutions, we desire to offset the spheres. From many points of view, the standard bubbles and Delaunay-Fowler-type solutions exhibit remarkable rigidity. Nevertheless, they can be moved by conformal transformations.


\vspace{0.3in}

{\bf \large {\bf  3. \ \ Double Kelvin transform being the identity map}}

\vspace{0.2in}

Define the reflection on the sphere with center at $\xi$ and radius $a > 0$ by 
$$
R_{\xi, \,a} (x) := \xi + {{a^2 (x - \xi)}\over {|x - \xi|^2}} \mfor x \not= \xi\,. \leqno (3.1)
$$
It is direct to check (cf. (3.3) below) that $R_{\xi, \,a} (R_{\xi, \,a} (x)) = x$ for $x \not = \xi\,.$ We observe that  
$$
\xi + {{a^2 (x - \xi)}\over {|x - \xi|^2}} \not= \xi \ \ \ \ {\mbox{when}} \ \ x \not= \xi\,.
$$

\vspace*{0.1in}

\hspace*{0.5in}The Kelvin transform with center at $\xi$ and radius $a > 0$ of a function $u$ is given by 
$$
\tilde u_{\xi, \,a} (x) := {{a^{n - 2}}\over {|x - \xi|^{n - 2}}} \,\,u \left(\xi + {{a^2 (x - \xi)}\over {|x - \xi|^2}} \right) \mfor x \not= \xi\,. \leqno (3.2)
$$
For later applications, we remark that if $\xi$ and $a$ are fixed, and if $u (y) \gg 1$, then from (3.2), $\tilde u_{\xi, \,a} (x) \gg 1$ as well, where 
$y = R_{\xi, \,a} (x)\,.$ 
In addition, $R_{\xi, \,a}$ sends a set of small diameter not too close to $\xi$ to a set of small diameter.\bk
We verify that double Kelvin transform is the identity map. Let $\tilde{\tilde u}$ be the Kelvin transform of ${\tilde u}_{\xi, \,a}$ with center at $\xi$ and radius $a\,.$ We have
\begin{eqnarray*}
(3.3) \ \ \ \ \ \ \ \ \ \ \tilde{\tilde u} (x) & = & {{a^{n - 2}}\over {|x - \xi|^{n - 2}}} \,\, \tilde u \left(\xi + {{a^2 (x - \xi)}\over {|x - \xi|^2}} \right)\\
& = & {{a^{n - 2}}\over {|x - \xi|^{n - 2}}} \cdot {{ a^{n - 2} }\over { a^{ 2 (n - 2)} / |x - \xi |^{n - 2} }} \,\, u \left(\xi + {{a^2 \left( {{a^2 ( x - \xi)} \over { |x - \xi|^2}} \right)
}\over {a^4/|x - \xi|^2}} \right) \ \ \ \ \ \ \ \ \ \ \ \  \\
& = & u (x) \ \ \mfor x \not= \xi\,.
\end{eqnarray*}

\vspace{0.3in}

{\bf \large {\bf  4. \ \ Kelvin transform of standard bubbles}}

\vspace{0.2in}

Let 
$$
u_\lambda (x) = \left( {\lambda\over {\lambda^2 + |x - \xi|^2}} \right)^{{n - 2}\over 2} \mfor x \in {\R}^n. \leqno (4.1)
$$
As a special case of (3.1), the Kelvin transform of $u_\lambda$ with center at the origin and radius $a > 0$ is given by 
$$
\tilde u_\lambda (x) := {{a^{n - 2}}\over {|x|^{n - 2}}} \,\,u_\lambda \left({{a^2 x}\over {|x|^2}}\right) \mfor x \not= 0\,. \leqno (4.2)
$$
It satisfies the equation 
$$
\Delta \tilde u_\lambda  + n (n - 2) \tilde u_\lambda^{{n + 2}\over {n - 2}} = 0 \ \ \ \ {\mbox{in}} \ \ {\R}^n \setminus \{ 0 \}\,.
$$
Moreover, the origin is a removable singularity 
and hence $\tilde u_\lambda$ is a standard bubble. A direct pointwise calculation shows that \cite{Leung.7}
$$
{\tilde u}_\lambda (x) = \left( {{ {{\overline \lambda}}}\over {{\overline \lambda}^2 + |x - {\overline \xi}|^2}} \right)^{{n - 2}\over 2}
\mfor x \in {\R}^n, \leqno (4.3)
$$ 
where
$$
{\overline \xi} = {{a^2 \,\xi} \over { \lambda^2 + |\xi|^2}} \ \ \ \ {\mbox{and}} \ \ \ \ {\overline \lambda} = {{a^2\lambda}\over { \lambda^2 + |\xi|^2}}\,. \leqno (4.4)
$$
The key feature here is that, given $\xi$ and $\lambda$, we can choose the positive number  $a$ so that 
$$
{\overline \lambda} = {{a^2}\over { \lambda^2 + |\xi|^2}} = 1\,. \leqno (4.5)
$$
With such a choice of $a$ we have $\overline \lambda = \lambda$ and  
$\overline \xi = \xi\,.$ We may say that the standard bubble is symmetric under the specific Kelvin transform. Furthermore, $\overline \xi$ comes from 
the reflection of a point $\xi'$ which is governed by the relation
$$
 {{a^2 \xi'}\over {|\xi'|^2}}  = \overline \xi = {{a^2 \xi}\over {\lambda^2 + |\xi|^2 }}\ \  \Longrightarrow \ \ |\xi'| = |\xi| + {{\lambda^2}\over {|\xi|}}\,. \leqno (4.6)
$$
(4.6) shows that the point $\xi'$ is offset from the center $\xi$ by the amount of $\lambda^2/|\xi|\,.$

\vspace{0.3in}

{\bf \large {\bf  5. \ \ Kelvin transform of the canonical standard bubble}}

\vspace{0.2in}

Consider the canonical standard bubble 
$$
u_s (x) = \left( {1\over {1 + |x|^2}} \right)^{{n - 2}\over 2} \mfor x \in {\R}^n.  $$
We look at the dual situation of the case in the previous section , that is, the center of the standard bubble is at the origin and the sphere for reflection has its center at a point $\xi \in {\R}^n.$ 
From (3.1), the Kelvin transforms of $u_s$ with center at $\xi$ and radius $a > 0$ is given by 
\begin{eqnarray*}
(5.1) \ \ \ \ \ \ \ \ \ \ \ \ \ \ \  \tilde u_{\xi, \,a} (x) & := & {{a^{n - 2}}\over {|x - \xi|^{n - 2}}} \,\,u_s \left(\xi + {{a^2 (x - \xi)}\over {|x - \xi|^2}} \right)\\
& = & {{a^{n - 2}}\over {|y|^{n - 2}}} \,\,u_s \left(\xi + {{a^2 y}\over {|y|^2}} \right) \mfor x \not= \xi\,,  \ \ \ \ \ \  \ \ \ \ \ \ \ \ \ \ \ \ \ \ \ \ \ \ \ \ 
\end{eqnarray*}
where $y = x - \xi \not= 0\,.$ The function 
$$
{{a^{n - 2}}\over {|y|^{n - 2}}} \,\,u_s \left(\xi + {{a^2 y}\over {|y|^2}} \right) \mfor y \not= 0 \leqno (5.2)
$$
is the same as the Kelvin transforms with center at the origin and radius $a > 0$ of the standard bubble    
$$
\bar u (y) :=  \left( {1\over {1 + |y - (- \xi)|^2}} \right)^{{n - 2}\over 2} \mfor y \in {\R}^n, 
\leqno (5.3)
$$
which is given by $\displaystyle{ \,{{a^{n - 2}}\over |y|^{n - 2}}\,\, \bar u \left( {{a^2 y}\over {|y|^2}} \right).\,}$ From section 4, we obtain  
\begin{eqnarray*}
(5.4) \ \ \ \ \ \ \ \ \ \ \ \ \ \ \ \tilde u_{\xi, \,a} (x)& =&\left( {{ {{\overline \lambda}}}\over {{\overline \lambda}^2 + |y - ({\overline{- \xi}}))|^2 }} \right)^{{n - 2}\over 2}\\
& = & \left( {{ {{\overline \lambda}}}\over {{\overline \lambda}^2 + \bigg\vert x - \left(\xi + ({\overline{- \xi}})\right)\bigg\vert^2}} \right)^{{n - 2}\over 2}\\
& = & \left( {{ {{\overline \lambda}}}\over {{\overline \lambda}^2 + |x -  \xi_c|^2}} \right)^{{n - 2}\over 2}\mfor x \not= \xi\,, \ \ \ \ \ \ \ \  \ \ \ \ \ \ \ \ \ \ \ \ \ \ \    
\end{eqnarray*}
where, by (4.4), 
$$
\xi_c : = \xi + ({\overline{-\xi}}) = \xi + {{a^2 \, (-\xi)  } \over { 1 + |\xi|^2}} \ \ \ \ {\mbox{and}} \ \ \ \ {\overline \lambda} = {{a^2}\over {1 + |\xi|^2}}\,. \leqno (5.5)
$$
 $\tilde u_{\xi, \,a}$ has a removable singularity at $\xi$ and is recognized as the standard bubble with center at $\xi_c\,.$ Choose the positive number $a$ so that $\overline \lambda = 1$, that is, 
$$
{a^2\over{1 + |\xi|^2}} = 1\,.
$$ 
It follows that $\xi_c = \xi - \xi = 0$. With such a choice of $a\,,$  $\tilde u_{\xi, \,a} = u_s$ in ${\R}^n$, which is the symmetry we discuss in section 4. The origin of $\tilde u_{\xi, \,a}$ comes from a point $x_c$ via the reflection which can be found by 
putting $x = \xi_c = 0$ in the expression 
$$
\xi + {{a^2 (x - \xi)}\over {|x - \xi|^2}}\,. \leqno (5.6)
$$
This gives us
$$
x_c = \xi + {{a^2 (0 -\xi)}\over {|0 - \xi|^2}} = \xi - {{(1 + |\xi|^2)\, \xi}\over {|\xi|^2}} = - {{\xi}\over {|\xi|^2}}\,, \leqno (5.7)
$$
where  (5.6) is used. Thus the {\it center} of the transformed standard bubble $\tilde u_{\xi, \,a}$ in (5.4)  is the Kelvin transform of a point which is offset from the origin by a distance of 
$1/|\xi| \,$ (cf. (4.6)).\bk
Another way to see this is to  solve for the point $x_c$ on the line segment in the direction of $\xi$ from the equation 
$$ 
\xi + {{a^2 (x_c - \xi)}\over {|x_c - \xi|^2}} =  \xi_c = 0
$$
We obtain
$$
 |x_c - \xi| = {{a^2}\over {|\xi|}} = {{1 + |\xi|^2}\over {|\xi|}} 
\Longrightarrow |x_c - \xi|  = |\xi| + {1\over {|\xi|}} \Longrightarrow x_c = - {\xi\over {|\xi|^2}}\,.  \leqno (5.8)
$$
As the application of the same reflection twice gives the identity map (except at the center of the sphere), the methods above, as expected, lead to the same result.

\newpage

\hspace*{0.5in}The reflection $R_{\xi, a}$ defined in (3.1) sends an open ball to an open ball. Denote by $B_o (\delta)$ the open ball with center at the origin and radius $\delta$ ($\delta > 0$). Then $R_{\xi, a} (B_o (\delta))$ is an open ball with radius $\breve \delta$ given by 
$$
\breve \delta = {1\over 2} \left[ {{a^2}\over {|\xi| - \delta}} - {{a^2}\over {|\xi| + \delta}} \right] = {{a^2 \delta}\over { |\xi|^2 - \delta^2}}\,. \leqno (5.9)
$$
Taking into the account that $a^2 = 1 + |\xi|^2$, we obtain 
$$
\breve \delta = {{(1 + |\xi|^2)\,\delta}\over { |\xi|^2 - \delta^2}} = \delta \cdot \left( {{1 + |\xi|^{-2} }\over{1 - \delta^2/|\xi|^2}} \right)\,. \leqno (5.10)
$$
In particular, if $|\xi| \gg 1$ and $\delta$ is small, then $\breve \delta \approx \delta\,.$ Indeed, when $|\xi| \gg 1\,,$ 
$$
a = \sqrt{1 + |\xi|^2\,} \approx |\xi| + {1\over {2 |\xi|}}\,. \leqno (5.11)
$$
The effect of the reflection $R_{\xi, a}$ on a small neighborhood of the origin can be approximated by the reflection upon the hyperplane that is perpendicular to $\xi$ and has distance $a$ from $\xi$. One can then link up (albeit only roughly) (5.7), (5.10) and (5.11) through the  geometric  picture.


\vspace{0.3in}

{\bf \large {\bf  6. \ \ Construction}}

\vspace{0.2in}

We begin with the basic procedure to offset the center of the canonical standard bubble and insert a modified Delaunay-Fowler-type   solution on a small open set, and repeat the construction infinitely many times. The result is a function whose graph can be visualized as a lot of long and thin needles stuck uprightly on a piece of cotton, getting very `dense' near the origin.\bk 
Start with the canonical standard bubble $u_s$ and a point $\xi \in {\R}^n$ with $|\xi| \gg 1$. We apply 
the Kelvin transforms  with center at $\xi$ and radius $a > 0$ to $u_s$ (see (5.1)), and choose $a$ such that $a^2 = 1 + |\xi|^2\,.$ 
The transformed standard bubble has a removable singularity at $\xi\,.$ From section 5, 
$\,\,\tilde u_{\xi, \,a} = u_s$ in ${\R}^n.$ 
The center of $\tilde u_{\xi, \,a}$ is the reflection via (3.1)  of a point $x_c$ with $|x_c| = 1/|\xi|$ ((5.7) and (5.8)). By choosing $\xi$ with $|\xi| \gg 1$, $x_c$ is close to the origin.\bk
We may modify $\tilde u_{\xi, \,a}$ on a small neighborhood of the origin as described in section 2 and obtain a smooth positive function $\breve u_{\xi, \,a}$ so that 
$$
\breve u_{\xi, \,a} (x) = \tilde u_{\xi, \,a} (x)  \mfor |x| \ge e^{-D}\,, \leqno (6.1)
$$
and $\breve u_{\xi, \,a}$ satisfies that equation 
$$
\Delta \breve u_{\xi, \,a} + n (n - 2) \breve K {\breve u}_{\xi, \,a}^{{n + 2}\over {n - 2}} = 0 \ \ \ \ {\mbox{in}} \ \ {\R}^n\,,  \leqno (6.2)
$$
where $\breve K$ is a smooth function on ${\R}^n$ with $|\breve K - 1|\le \varepsilon$ in ${\R}^n.$ Here $D$ and $\varepsilon$ are arbitrarily given positive numbers. Furthermore,  $\breve u_{\xi, \,a} (0)$ can be made as large as we like, either by taking a large $T$ or by incorporating more cycles.\bk
The Kelvin transform of $\breve u_{\xi, \,a}$ with center at $\xi$ and radius $a > 0$, by (3.3), gives us a function ${\breve{\breve u}}_\xi$ with the following properties:
$$
{\breve{\breve u}}_\xi (x) = u_s (x) 
$$
for $x$ outside a small open neighborhood $U_1$ of $x_c$ (by choosing $D$ large enough, say, $e^{D} \gg |\xi|$, we may assume that $0 \not\in {\overline U}_1$; see also (5.10)). The measure of $U_1$ can be made arbitrarily small by selecting $D \gg 1$ and $|\xi| \gg 1$. Furthermore, by taking $\breve u_{\xi, \,a} (0)$ to be large (cf. the remark following (3.2) in section 3), we have 
$$
{\breve{\breve u}}_\xi (x_c) \gg |x_c|^{{2-n}\over 2}\,.  \leqno (6.3)
$$ 
As ${\breve{\breve u}}_\xi = u_s$ outside $U_1\,,$ ${\breve{\breve u}}_\xi$ has a removable singularity at $\xi$ and satisfies the equation 
$$
\Delta {\breve{\breve u}}_\xi + n (n - 2) \,\breve{\breve K} \, {\breve{\breve u}}_\xi^{{n + 2}\over {n - 2}} = 0\,,  \leqno (6.4)
$$
where $\breve{\breve K}$ is a smooth function on ${\R}^n$ with 
$$
\breve{\breve K} (x) := \breve K \left(\xi + {{a^2 (x - \xi)}\over {|x - \xi|^2}} \right) \mfor x \not= \xi\,.  \leqno (6.5)
$$ 
Hence $|\breve{\breve K} - 1|\le \varepsilon$ in ${\R}^n.$ Moreover, $\breve{\breve K} \equiv 1$ outside $U_1\,.$ \bk
The second step is to repeat the above construction by taking $x_c$ closer to the origin, and the open neighborhood $U_2$  around $x_c$ to be smaller so that $U_2 \cap U_1 = \emptyset$ and $0 \not\in {\overline U}_2\,.$   
Here $U_2$ carries the meaning of $U_1$ in the procedure. 
As there is no modification occuring inside $U_1$ in this step, after the double Kelvin transform, nothing is changed inside $U_1\,.$\bk
Continue the above construction by choosing larger and larger $|\xi|$ so as to move the offset centers closer and closer to the origin, larger $D$ to make the neighborhood of the offset center $x_c$ smaller, and larger $T$ to enlarge the value of the modified function at $x_c$, and to keep check of the inequality $|\breve{\breve K} - 1|\le \varepsilon$ in ${\R}^n.$
Hence we can define a smooth positive function $u_b$ on ${\R}^n \setminus \{ 0 \}$ with the following properties. There exists a sequence $\{ x_i \}$ with $|x_i| \to 0$ as $i \to \infty$ such that 
$u_b (x_i) \to \infty$ in an arbitrary fast rate (for example, $\displaystyle{u_b (x_i) \,|x_i|^{(n-2)/2} \to \infty}\,$ as $\,i \to \infty$). Furthermore, 
$u_b$ satisfies the equation 
$$
\Delta u_b + n (n - 2) \,K u_b^{{n + 2}\over {n - 2}} = 0 \ \ \ \ {\mbox{in}} \ \ {\R}^n \setminus \{ 0 \}\,, \leqno (6.6)
$$ 
where $K$ is a smooth function on ${\R}^n \setminus \{ 0 \}$ with $|K - 1| \le \varepsilon$ in ${\R}^n \setminus \{ 0 \}\,.$ Furthermore, $K \equiv 1$ and $u_b \equiv u_s$ outside a small open set $U$ which is the union of the above neighborhoods around the offset centers. The  measure of $U$ can be made arbitrarily small.

\vspace{0.3in}

{\bf \large {\bf  7. \ \ The Lipschitz condition}}

\vspace{0.2in}

When $n > 4\,,$ we show how to select the parameters $T$, $D$ and $\xi$ suitably so that $K$ in (6.6) can be extended as a Lipschitz function on ${\R}^n.$ 
To better keep track of the change on $K\,,$  we limit ourselves only to two cycles in the modification described in section 2. Following the notations in that section, we consider the function $\breve v_T$ with 
$$
{\breve{v}}_T (t) = v_s (t) \mfor t \le  D\,,  \leqno (7.1)
$$
$$
{\breve{v}}_T (t) = v_s (t - T) \mfor t \ge T - D\,, \leqno (7.2)
$$
and
$$
{\breve{v}}_T'' - {{(n - 2)^2}\over 4}\, {\breve{v}}_T + n (n - 2) {\breve K}\, {\breve{v}}_T^{{n + 2}\over {n - 2}} = 0 \ \ \ \ {\mbox{in}} \ \ \R\,. \leqno (7.3)
$$
Here $T \gg D\,.$ 
Accordingly, 
$$
\breve K (t) = 1  \mfor t \le D\,, \ \ t \ge T - D\,.
$$ 
Let 
$\breve u_T$ be the transform via (2.13) of $\breve v_T\,,$ and
$$
K (x) := {\breve K} (t, \theta)\,, \ \ \ \ {\mbox{where}} \ \ \  e^{-t} = x \ \ \ {\mbox{and}} \ \ \theta = {x\over {|x|}} \mfor x \not= 0\,.
$$
(For clarity sake we still use $K$ as the notation here, but mindful of (6.6).)
We have
$$
\Delta \breve u_T + n (n - 2) K \,{\breve u}_T^{{n + 2}\over {n - 2}} = 0 \ \ \ \ {\mbox{in}} \ \ {\R}^n \setminus \{ 0 \}\,, \leqno (7.4)
$$
and
$$
K (x) = 1 \mfor |x| \ge e^{-D}, \ \ \ \ 0 <  |x| \le e^{-T + D}.
$$ 
Define $K (0) = 1\,.$ We seek to show that, for $T$ large enough {\it and\,} $n > 4\,,$ 
$$
|K (x) - K (y)| \le |x - y| \mfor {\mbox{all}} \ \ x\,, \ y \in {\R}^n. \leqno (7.5)
$$
In case $y = 0\,,$ then (7.5) clearly holds for $|x| \ge e^{-D}$, $|x| \le e^{-T + D}\,.$  Consider the case when 
$$
y = 0 \ \ {\mbox{and}} \ \ e^{-T + D} \le |x| \le e^{-D}\,.
$$ 
By proposition 2.4, (2.8), and the construction, we have
$$
|K (x) - 1| \,\,\le\,\, C_1 (D) \,\eta_T^2 \,\, \le \,\,C_2 (D) \,e^{-  \left( {{n - 2}\over 2} \right) T}\,. \leqno (7.6)
$$
Here $C_1 (D)$ and $C_2 (D)$ are positive constants that depend on $D$ and $n$ only, with emphasis on $D$. 
Set $\displaystyle{c = {{n - 4}\over 2}}$ so that $c > 0$ when $n > 4\,.$ From (7.6) we obtain 
$$
|K (x) - 1| \le  \left[ {{C_2 (D)}\over {e^D}} \right] e^{-(1 + c) T + D} = \left[ {{C_2 (D)}\over {e^D\, e^{cT}}} \right] e^{- T + D} 
\le \left[ {{C_2 (D)}\over {e^{c T + D}}}  \right]\, |x|\,.  \leqno (7.7)
$$ 
Choosing $T$ to be large, we have $|K (x) - 1| \le |x|\,$ for $\,e^{-T + D} \le |x| \le e^{-D}.$ Hence (7.5) holds for $y = 0$ and $x\in {\R}^n.$\bk
Let $x, \ y \in {\R}^n \setminus \{ 0 \}\,.$ If $|x| < e^{-T + D}$ and $|y| < e^{-T + D}\,,$ 
then 
$$
|K (x) - K (y)| = |1 - 1| = 0 \le |x - y|\,.
$$ 
Similarly for the subcases (i) $|x|>  e^{-D}$ and $|y| >  e^{-D};$  (ii) $|x|>  e^{-D}$ and $|y| < e^{-T + D};$ and (iii) $|x| < e^{-T + D}$ and $|y| >  e^{-D}.$  We need only to consider the case when  $e^{-T + D} \le |x| \le e^{-D}$ (the case $e^{-T + D} \le |y| \le e^{-D}$ is  similar). Let $s$ and $t$ be positive numbers such that 
$|x| = e^{-s}$ and $|y| = e^{-t}\,.$ Thus $D \le s \le T - D\,.$

\newpage

\hspace*{0.5in}Assume that $t \ge s\,.$ As $K$ is rotationally symmetric about the origin, we may, if necessary,  replace $y$ by another point $y'$ with $|y| = |y'| = e^{-t}$ and $y'/|y'| = x/|x|\,.$ Furthermore,  
$$
|K (x) - K (y)| = |K (x) - K (y')| \ \ \ \ {\mbox{and}} \ \ |x - y| \ge |x - y'|\,.
$$
Thus if $|K (x) - K (y')| \le |x - y'|$, then we also have $|K (x) - K (y)| \le |x - y|\,.$ 
Hence we may assume that $y/|y| = x/|x|\,,$ and seek to establish (7.5).\bk
As ${\breve K} = 1$ on $(T- D\,, \ \infty)\,,$ we may make $t$ smaller until $t \le T- D\,,$ if necessary, without changing the left hand side of (7.5).   By doing so, the right hand side of (7.5) actually decreases. Hence we may also assume that 
$t \le T - D\,.$\bk
Using the triangle inequality again and the inequality $\displaystyle{e^{t - s} \ge 1 + (t - s)}$ for $t - s \ge 0\,,$ we verify that  
$$
|x - y| \ge  e^{-s} - e^{-t} \ge e^{-t} \left[e^{(t - s)} -1 \right] \ge e^{-T + D} (t - s)\,,
\leqno (7.8)
$$
as $\,t \le T - D\,.$ 
We have 
\begin{eqnarray*}
(7.9) \ \ \ \ \ \ \ \ \ \ \ 
|K (x) - K (y)| & = & |\breve K (s) - \breve K (t)| \le C_3 (D)\, e^{-(1 + c) T} |t - s|
\ \ \ \ \ \ \ \ \ \ \ \ \ \ \ \ \ \ \\
& = & {{C_3 (D)}\over {e^{cT+D}}} \, e^{-T + D} (t - s) \le |x - y|  
\end{eqnarray*}
for $T \gg 1\,,$ where we use proposition 2.4, (2.8) and (7.8). Here $C_3 (D)$  is a positive constant that depends on $D$ and $n$ only.\bk
As for the case when $s > t$ and $D \le s \le T - D\,,$ we need only to change the positions of $s$ and $t$ in (7.8) and (7.9). Likewise, we treat the case $D \le t \le T - D\,.$\bk
Take the Kelvin transform with center at $\xi$ and radius $a$ on $\breve u_T$, with $a^2 = 1 + |\xi|^2$ and $|\xi| \gg 1$. As $K \not= 1$ only in a small neighborhood of the origin, and $|\xi| \gg 1$ means that the reflection is close the reflection of a hyperplane (cf. the remark following (5.10)), we see that the Kelvin transform does not change the essential form of the Lipschitz condition on $K$ in (7.5).\bk 
We apply the above argument to the remaining constructions as described in section 6 (each step limited to two cycles only). In the final stage the function $K$ in (6.6),\, after taking into the account of the triangle inequality and choosing all the $T,$ $\xi$ and $D$ suitably (cf. (7.7) and (7.9)), satisfies
$$ 
|K (x) - K (y)| \le C\,|x - y| \mfor x, \ y \in {\R}^n \ \ \ \ (n > 4\,; \ \ K (0) := 1)\,.  \leqno (7.10)
$$
Here $C$ is a positive constant. Given a number $\alpha \in (0, 1]$, we have $|x - y|^\alpha \ge |x - y|$ for $|x - y| \le 1\,$, and 
$$
4 \ge |x - y| \ge 1 \ \ \Longrightarrow |x - y| \le (|x - y|^{1-\alpha} )|x - y|^\alpha \le 4\,|x - y|^\alpha\,. 
$$
Hence we obtain 
$$
|K (x) - K (y)| \le 4\, C |x - y|^\alpha \mfor x, \ y \in {\bar B}_2\,.
$$

\vspace*{0.2in}

{\bf Remark 7.11.} \ \  It is unlikely that the term $e^D$, which can be found in the denominators in (7.7) and (7.9), is large enough to cancel the dependence of $D$ in $C_2 (D)$ and $C_3 (D)\,.$ This can be seen in (A. 10) below, where the constant $C_1 (D)$ is in the order of $e^{(n + 2)D}$ (cf. also (2.6)). Therefore the presence of $e^{cT}$ in the (7.7) and (7.9) is rather necessary, i.e., $n$ is required to be bigger than $4\,.$


\vspace{0.3in}

{\bf \large {\bf Appendix}}

\vspace{0.2in}

The purpose is to derive estimate (2.8). Keeping the notations used in section 2, we begin  with the inequalities 
$$
| v_T (t) - v_s (t)|  \,\le\, C_D \,\eta_T^2 \ \ \ \ {\mbox{and}} \ \ \ \  
| v'_T (t) - v'_s (t)|  \,\le\, C_D \,\eta_T^2 \mfor t \in [-2D, \ 2D] \leqno (A.1)
$$
and for all $T$ large enough. Here $C_D$ is a positive constant that depends on $D$ and $n$ only. Similar estimates are shown in \cite{Leung.8} on a small interval for practical reasons, and they can be extended by a continuation argument (cf. inequality (2.18) in \cite{Leung.8}).\bk 
Given positive numbers $\alpha\,,$ $\delta$ and $\hat c\,,$ where $\delta \le \hat c\,,$ the inequality 
$$
|\,x^\alpha - y^\alpha\,| \le C\,|x - y| \ \ \ \ {\mbox{holds \ \ for}} \ \ \  \delta \le x \le \hat c, \ \, \delta \le  y \le \hat c\,, \leqno (A.2)
$$
with a positive constant $C$ that depends on $\alpha\,,$ $\hat c$ and $\delta$ only. (A.2) can be established by an integration method (cf. \cite{Leung.8}).  
The Delaunay-Fowler-type   solution $v_T$ is uniformly bounded from above and away from zero in $[-2D, \ 2D]$ 
for all $T \gg 1\,$ (see \cite{M-P} and \cite{Leung.8}).\bk
Let $\phi_1$ be a smooth function on $\R$ such that $0 \le \phi \le 1$ in $\R$ and 
\[ \phi_1 (t) \ = \ \left\{ \begin{array}   
{r@{\quad \mbox{for}\quad}l}
 1 & \ \  t \le D\,,\\  
0  & \ \  t \ge 2D\,. 
\end{array} \right.  \]
Take $D \ge 1\,.$ We also require that 
$$
|\,\phi_1' (t)| \,\le 2 D\,, \ \ 
|\,\phi_1'' (t)| \,\le \,2 D\,, \ \ \ \ {\mbox{and}} \ \ \  |\,\phi_1''' (t)| \,\le \,2 D  \mfor t \in (D\,, \,2D)\,. \leqno (A.3)
$$
Set $\phi_2 := 1 - \phi_1$ in $\R$. Define
$$
{\breve v}_T := \phi_1 \,v_s + \phi_2 \,v_T\ \ \ \ {\mbox{in}} \ \ \R\,, \leqno (A.4)
$$
where $T \gg 1\,.$ We have ${\breve v}_T (t) = v_s (t)$ for $t \le D$ and 
${\breve v}_T (t) = v_T (t)$ for $t \ge 2D\,.$ 
${\breve v}_T$ satisfies the equation
$$
{\breve v}_T'' -  {{(n - 2)^2}\over 4} \,{\breve v}_T (t) + n (n - 2) {\breve K} {\breve v}_T^{{n + 2}\over {n - 2}} (t) = 0 \ \ \ \ {\mbox{in}} \ \ \R\,, \leqno (A.5)
$$
where
\begin{eqnarray*}
(A.6)\ \,\, 
{\breve K} (t) & := & \left[- {\breve v}_T'' (t)  + {{(n - 2)^2}\over 4} \,{\breve v}_T (t)\right]
\,\left[ n (n - 2)\, {\breve v}_T^{{n + 2}\over {n - 2}} (t) \right]^{-1}\\
& = & \{
n (n - 2)  \left[ \phi_1 (t) \, v_s^{{n + 2}\over {n - 2}} (t) \,+ \,
\phi_2 (t) \,v_T^{{n + 2}\over {n - 2}} (t) \right] + \phi_1' (t) \,[v_T' (t) 
- v_s' (t) ]\\
& \ & \ \ \ \ \ \ \ \ \ \ \ \ \ \ \ \ \ \ \ \
 \,\,+ \,\, \phi_1'' (t) \,[v_T (t) 
- v_s (t) ] \}
\,\left[ n (n - 2)\, {\breve v}_T^{{n + 2}\over {n - 2}} (t) \right]^{-1} \ \ \ \ \ \ 
\end{eqnarray*}
for $t \in \R\,.$ Here we use (A.4), (A.5), and equation (2.2), which $v_s$ and $v_T$ satisfy.
It follows that 
\begin{eqnarray*}
(A.7) & \, &{\breve K}' (t)  =  \{
n (n - 2) \left[  {{n + 2}\over {n - 2}} \left( \phi_1 (t) \,v_s^{4\over {n - 2}} (t) \,v_s' (t) \,+ \,
\phi_2 (t) \,v_T^{4\over {n - 2}} (t)\, v_T' (t) \right) \right. \ \ \ \  \\
& \ & \ \ \ \ \ \ \ \ \ \ \ \ \ \left. +\,\, \phi_1'(t) \,v_s^{{n + 2}\over {n - 2}} (t) \,+ \,
\phi_2' (t)\,v_T^{{n + 2}\over {n - 2}} (t) \right]\\
& \ &  \ \ \ \ \ \ + \,\,2 \phi_1'' (t) \,[v_T' (t) 
- v_s' (t) ]  + \phi_1' (t) \,[v_T'' (t) 
- v_s'' (t) ]\\
& \ &  \ \ \ \ \ \ + \,\, \phi_1''' (t) \,[v_T (t) 
- v_s (t) ] \,\}
\,\left[ n (n - 2)\, {\breve v}_T^{{n + 2}\over {n - 2}} (t) \right]^{-1}\\
& \ &  \ \ \ \ \ \  -\, \{
n (n - 2)  \left[ \phi_1 (t) \,v_s^{{n + 2}\over {n - 2}} (t) \,+ \,
\phi_2 (t) \,v_T^{{n + 2}\over {n - 2}} (t) \right] + \phi_1' (t) \,[v_T' (t) 
- v_s' (t) ]  \\
& \ &  \ \ \ \ \ \ +  \,\, \phi_1'' (t) \,[v_T (t) 
- v_s (t) ] \}
\, \left[ {{n + 2}\over {n (n - 2)^2}} \right] \,
\left[  {\breve v}_T^{{n + 2}\over {n - 2}} (t) \right]^{-2} \left[   {\breve v}_T^{4\over {n - 2}} (t)\, {\breve v}_T' (t) \right] \ \ \ \ \ \ \ \ 
\end{eqnarray*}
for $t \in \R\,.$
Equation (2.2) shows that 
$$
1 = \left[- v_s'' (t)  + {{(n - 2)^2}\over 4} \,v_s (t)\right]
\,\left[ n (n - 2)\, v_s^{{n + 2}\over {n - 2}} (t) \right]^{-1} \mfor t \in \R\,. \leqno (A.8)
$$
After differentiating both sides of (A.8) we obtain
\begin{eqnarray*}
(A.9) & \ &\\
0 & = & \left[- v_s''' (t)  + {{(n - 2)^2}\over 4} \,v_s' (t)\right]
\,\left[ n (n - 2)\, v_s^{{n + 2}\over {n - 2}} (t) \right]^{-1}\\
& \ &  - \left[- v_s'' (t)  + {{(n - 2)^2}\over 4} \,v_s (t)\right]
\,\left[ {{n + 2}\over {n (n - 2)^2}} \right] \,
\left[  v_s^{{n + 2}\over {n - 2}} (t) \right]^{-2} \left[   v_s^{4\over {n - 2}} (t) \,v_s' (t) \right]\\
& = & \left[n (n - 2) \left( {{n + 2}\over {n - 2}} v_s^{4\over {n - 2}} (t) \,v_s' (t) \right) \right]
\,\left[ n (n - 2)\, v_s^{{n + 2}\over {n - 2}} (t) \right]^{-1}\\
& \ & \ \ \ \  - \left[n (n - 2)\,v_s^{{n + 2}\over {n - 2}} (t)\right]
\,\left[ {{n + 2}\over {n (n - 2)^2}} \right] \,
\left[  v_s^{{n + 2}\over {n - 2}} (t) \right]^{-2} \left[   v_s^{4\over {n - 2}} (t) \,v_s' (t) \right] \ \ \ \ \ \ \ \ \ \ \ 
\end{eqnarray*}
for $t \in \R\,.$ Here we use equation (2.2) to find $v_s'''.$ 
Subtracting (A.7) and (A.9) and consider the leading terms in (A.7) and (A.9) (ignoring the constant $(n + 2)/(n - 2)$ in front), we have
\begin{eqnarray*}
& \ &  \Bigg\vert 
{{\phi_1 (t) \,v_s^{4\over {n - 2}} (t) \,v_s' (t) \,+ \,
\phi_2 (t) \,v_T^{4\over {n - 2}} (t)\, v_T' (t) }\over {{\breve v}_T^{{n + 2}\over {n - 2}} (t)}} - {{ v_s^{4\over {n - 2}} (t) v_s' (t)}\over {v_s^{{n + 2}\over {n - 2}} (t)}} \Bigg\vert\\
& = &  \Bigg\vert \,
{{ v_s^{{n + 2}\over {n - 2}} (t)\, [\phi_1 (t)  \,v_s^{4\over {n - 2}} (t) \,v_s' (t) \,+ \,
\phi_2 (t) \,v_T^{4\over {n - 2}} (t)\, v' (t)] - {\breve v}_T^{{n + 2}\over {n - 2}} (t) \, v_s^{4\over {n - 2}} (t) v_s' (t)}\over { ({\breve v}_T\, v_s)^{{n + 2}\over {n - 2}} (t)}}\, \Bigg\vert\\
& \le & C_1 (D) \, \Bigg\vert \,
 v_s^{{n + 2}\over {n - 2}} (t)\, \left[\phi_1 (t)  \,v_s^{4\over {n - 2}} (t) \,v_s' (t) \,+ \,
\phi_2 (t) \,v_T^{4\over {n - 2}} (t)\, v_T' (t)\right]\\
& \ & \ \ \ \ \ \ \ \ \ \ \ \ \ \ \ \ \ \ \ \ \ \ \ \ \ \   - \,\,{\breve v}_T^{{n + 2}\over {n - 2}} (t) \, v_s^{4\over {n - 2}} (t)\, v_s' (t) \, \Bigg\vert\\
& = & C_1 (D) \, \Bigg\vert \,
 v_s^{{n + 2}\over {n - 2}} (t)\, \left[\phi_1 (t)  \,v_s^{4\over {n - 2}} (t) \,v_s' (t) \,+ \,
\phi_2 (t) \,v_T^{4\over {n - 2}} (t)\, v_T' (t) - v_s^{4\over {n - 2}} (t) v_s' (t)\right]\\
& \ & \ \ \ \ \ \ \ \ \ \ \ \ \ \ \ \ \ \ \ \ \ \ \ \ \ \   - \,\,[v_T^{{n + 2}\over {n - 2}} (t) -  v_s^{{n + 2}\over {n - 2}} (t)]\, v_s^{4\over {n - 2}} (t) \,v_s' (t) \, \Bigg\vert\\
& \le & C_1 (D) \, \left\{ \Bigg\vert \,
 v_s^{{n + 2}\over {n - 2}} (t)\, \left[ 
\phi_2 (t) \,v_T^{4\over {n - 2}} (t)\, v_T' (t) -\phi_2 (t) \,v_s^{4\over {n - 2}} (t) \,v_s' (t) \right] \Bigg\vert \right.\\
& \ & \ \ \ \ \ \ \ \ \ \ \ \ \ \  
 + \, \left.  \Bigg\vert v_T^{{n + 2}\over {n - 2}} (t) -  v_s^{{n + 2}\over {n - 2}} (t) \Bigg\vert \, \Bigg\vert v_s^{4\over {n - 2}} (t) v_s' (t) \, \Bigg\vert \right\}\\
& \le & C_2 (D) \left\{ \Bigg\vert v_T^{4\over {n - 2}} (t)\, v_T' (t) - \,v_s^{4\over {n - 2}} (t) \,v_s' (t) \Bigg\vert  +  \Bigg\vert v_T^{{n + 2}\over {n - 2}} (t) -  v_s^{{n + 2}\over {n - 2}} (t) \Bigg\vert \right\}\\
& \le & C_2 (D) \left\{   v_T^{4\over {n - 2}} (t) | v_T' (t) - v_s' (t)| + |v_s' (t)|\, \bigg\vert  v_T^{4\over {n - 2}} (t) - v_s^{4\over {n - 2}} (t) \bigg\vert \right.\\
& \ & \ \ \ \ \ \ \ \ \ \ \ \ \ \  +\, \left. 
\Bigg\vert v_T^{{n + 2}\over {n - 2}} (t) -  v_s^{{n + 2}\over {n - 2}} (t) \Bigg\vert  \right\}\\
& \le & C_3 (D) \, ( \,| v_T' (t) - v_s' (t)|  + | v_T (t) - v_s (t)|\,)  \le C_4 (D) \,\eta_T^2 \mfor t \in [D, \ 2D]\,.  
\end{eqnarray*} 
Here we use (A.1) and (A.2).  In addition,  $C_1$, $C_2,$ $C_3$ and $C_4$ are positive constants that depend on $D$ and $n$ only. We may take, for examples, 
$$
C_1 (D) = 
\max_{[D, \ 2D]} \, \left\{ \,{1\over {
({\breve v}_T \,v_s)^{{n + 2}\over {n - 2}} }} \,\right\}\,, \leqno (A.10)
$$
and 
$$
C_2 (D) = \left(\, \max_{[D, \ 2D]} \left\{ v_s^{{n + 2}\over {n - 2}}, \ \,\,
v_s^{4\over {n - 2}}\,|v_s'|\,   \right\}\, \right) \cdot C_1 (D)\,.
$$
Similarly, we may define the value for 
 $C_3$ after 
taking into consideration of (A.2) and 
whether $4/(n - 2) \ge 1$ or not. 
Likewise, we consider the absolute values of the subtraction of the last third term in (A.7) with the last term in (A.8) (ignoring the constant $(n + 2)/(n - 2)$ again):
\begin{eqnarray*}
& \ &  \Bigg\vert \left[ \phi_1 (t)  \,v_s^{{n + 2}\over {n - 2}} (t) \,+ \,
\phi_2 (t) \,v_T^{{n + 2}\over {n - 2}} (t) \right] \left[ {{ {\breve v}_T' (t)}\over {{\breve v}_T^{{2n}\over {n - 2}} (t)}} \right] - \left[ \,v_s^{{n + 2}\over {n - 2}} (t)  \right] \left[ {{ v_s' (t)}\over {v_s^{{2n}\over {n - 2}} (t)}} \right] \Bigg\vert\\
& = &  \Bigg\vert {{ \left[ \phi_1 (t)  \,v_s^{{n + 2}\over {n - 2}} (t) \,+ \,
\phi_2 (t) \,v_T^{{n + 2}\over {n - 2}} (t) \right] \left[ v_s^{{2n}\over {n - 2}} (t) \right] [{\breve v}_T' (t)] - \left[ \,v_s^{{n + 2}\over {n - 2}} (t)  \right] \left[ {\breve v}_T^{{2n}\over {n - 2}} (t) \right] [v_s' (t)]}\over { ({\breve v}_T \,v_s)^{{2n}\over {n - 2}} (t) }} \Bigg\vert\\
& \le & C_5 (D) \left\{ \Bigg\vert \phi_1 (t)  \,v_s^{{n + 2}\over {n - 2}} (t) \,+ \,
\phi_2 (t) \,v_T^{{n + 2}\over {n - 2}} (t) - v_s^{{n + 2}\over {n - 2}} (t) \Bigg\vert \right.\\
& \ & \ \ \ \   \ \ \ \  \ \ \ \  \ \ \ \  \ \ \ \  \ \ \ \   \ \ \ \  \ \ \ \ \ \ \ \   \ \ \ \  \ \ \ \ \left. 
+ \, \,\Bigg\vert v_s^{{2n}\over {n - 2}} (t) - {\breve v}_T^{{2n}\over {n - 2}} (t) \Bigg\vert 
+ |  {\breve v}_T' - v_s' (t)| \right\}\\
& \le & C_6 (D) \,\eta_T^2 \ \ \mfor t \in [D, \ 2D]\,.  
\end{eqnarray*}
Observe in (A.7) that $\phi_1' = - \phi_2'.$ Hence we have 
\begin{eqnarray*}
& \ & |2 \phi_1'' (t)|\, \Bigg\vert \phi_1'(t) \,v_s^{{n + 2}\over {n - 2}} (t) \,+ \,
\phi_2'(t) \,v_T^{{n + 2}\over {n - 2}} (t) \Bigg\vert 
\left[ n (n - 2)\, {\breve v}_T^{{n + 2}\over {n - 2}} (t) \right]^{-1}\\
& \le & C_7 (D) 
\Bigg\vert  \,v_s^{{n + 2}\over {n - 2}} (t) \,- \,
 \,v_T^{{n + 2}\over {n - 2}} (t) \Bigg\vert \le C_8 (D) \,\eta_T^2 \mfor t \in [D, \ 2D]\,,
\end{eqnarray*}
where we use (A.1) and (A.2). 
It follows from equation (2.2) that 
\begin{eqnarray*}
& \ & |\phi_1' (t)| \,\bigg\vert v_T'' (t) - v_s'' (t) \bigg\vert \left[ n (n - 2)\, {\breve v}_T^{{n + 2}\over {n - 2}} (t) \right]^{-1}\\
& \le & C_9 (D) \Bigg\vert  {{(n - 2)^2}\over 4} v_T  - n (n - 2) v_T^{{n + 2}\over {n - 2}} (t)  - {{(n - 2)^2}\over 4} v_s + n (n - 2) v_s^{{n + 2}\over {n - 2}} (t) \Bigg\vert\\
& \le & C_{10} (D) \left\{ |v_T (t) - v_s (t)| + \Bigg\vert v_s^{{n + 2}\over {n - 2}} (t) - v_T^{{n + 2}\over {n - 2}} (t) \Bigg\vert \right\}\\
& \le & C_{11} (D)  |v_T (t) - v_s (t)| \le C_{12} (D) \,\eta_T^2 \mfor t \in [D, \ 2D]\,.
\end{eqnarray*}
As above, 
$$
|\phi_1''' (t)| \bigg\vert v_T (t) - v_s (t) \bigg\vert \left[ n (n - 2)\, {\breve v}_T^{{n + 2}\over {n - 2}} (t) \right]^{-1} \le C_{13} (D)  |v_T (t) - v_s (t)| \le C_{14} (D)\, \eta_T^2  
$$
for $t \in [D, \ 2D]\,.$  Akin to the above, we treat the fourth and the last terms  in (A.7). 
Hence we see that 
$$
|{\breve K}' (t)| = |{\breve K}' (t) - 0| \le C (D) \, \eta_T^2 \mfor t \in [D, \ 2D]\,. \leqno (A.11)
$$ 
Whilst (A.11) clearly holds outside $[D, \ 2D]\,.$ Upon integration we obtain
$$
|{\breve K} (t) - {\breve K} (s)| \le C (D) \, \eta_T^2 \, |t - s| \mfor s, \ t \in \R\,. \leqno (A.12)
$$
Here $C (D)$ and $C_i (D)$ are positive constants that depend on $D$ and $n$ only (for $T \gg 1$), and $i = 1, \ 2,\cdot \cdot \cdot , 14\,.$ Similar consideration can be applied on the cut-and-glue-in process on $(T- 2D,\  T - D)$ and elsewhere.

\pagebreak

\vspace*{1in}

{\bf Added.} \ \ After the paper was communicated to Steven Taliaferro, he informed the author that together with Lei Zhang they had just completed a paper which contained the following result.  Let $0 < b < (n-2)/2$  be a constant ($n \ge 3$). Then there exist positive $C^2$-functions $u$ which satisfy  
$$
(1-|x|^b)\,u^{{n + 2}\over {n - 2}} (x) \,\le \,-\Delta u (x) \,\le \,u^{{n + 2}\over {n - 2}} (x)  \ \ \ \ \mfor   0 < |x| <1\,, 
$$
and are arbitrarily large near the origin. \bk 
Here $(n - 2)/2$ appears to be a critical order of how fast  the term in the left hand side appraoches one.  It provides potentially fruitful insights into the study of equation (1.1).\bk 
Keeping the notations in section 7, it follows  from (7.7) that, for any {\it small}  positive number $\beta\,,$

\newpage

\begin{eqnarray*}
{\mbox{(a.1)}} \ \ \ \ \ \  \ \ \ \ \ \ \ \ \ |K (x) - 1| & \le & \left[ {{C_2 (D)}\over {e^{\left({{n - 2}\over 2} - \beta \right) D}}} \right] e^{-{{n - 2}\over 2} T + \left({{n - 2}\over 2} - \beta \right) D}\\
 & = &   \left[ {{C_2 (D)}\over {e^{\left({{n - 2}\over 2} - \beta \right) D} e^{\beta T}}} \right] e^{\left(-{{n - 2}\over 2} + \beta \right) T  + \left({{n - 2}\over 2} - \beta \right)D}
 \ \ \ \ \ \ \ \ \ \ \ \  \ \ \  \\
& \le & |x|^{ {{n - 2}\over 2} - \beta }
\end{eqnarray*}
for all  $T$ large enough (recall that $e^{-T + D} \le |x|$)\,.\, Furthermore, we need only $n \ge 3$ to draw the above conclusion. Thus the critical order $(n - 2)/2$ can also be found in our construction. The author thanks Taliaferro for informing him the interesting result, which led to the  consideration in (a.1).\bk 
Finally,  a natural question arises: Given a function  
$$
K \in C^1 (B_1 \setminus \{ 0 \}) \cap C^o (B_1)
$$
with 
$$
|K (x) - 1| \le C |x|^{{n - 2}\over 2} \mfor 0 < |x| < 1\,,
$$
where $C$ is a positive constant, does every positive $C^2$-solution $u$ of the equation 
$$
\Delta u + n (n - 2) \,K u^{{n + 2}\over {n - 2}} = 0 \ \ \ \ {\mbox{in}} \ \ B_1 \setminus \{ 0 \}
$$
have slow decay (1.2) near the origin? How about adding the condition that $\btd K \not = 0$ near the origin? (cf. \cite{Chen-Lin.2} and \cite{Lin.1}.)\bk
It is inspiring that the critical order $(n - 2)/2$ is related to the convergence rate  (2.8)  through proposition 2.4 on neck-sizes and periods of Delaunay-Fowler-type solutions. The same critical order is present in the constructions of exotic solutions in \cite{Chen-Lin.4} and \cite{Leung.8}, where the solutions unexpectedly breach a natural lower bound mirroring the slow decay estimate (1.2).


\end{document}